\documentclass[12pt]{article}
\usepackage{graphicx}
\usepackage{amsmath, amsthm, amssymb, amscd}
\usepackage{setspace}
\usepackage{mathrsfs}

\newtheorem{lemma}{Lemma}
\newtheorem{theorem}{Theorem}
\newtheorem*{corollary 1}{Corollary 1}

\title{\textbf{A new family of sharp conformally invariant integral inequalities}\thanks{The author would like to express his deep thanks to his supervisor Professor Robert McCann for a careful reading of the draft and many very useful discussions and suggestions. Without him, the work would not be possible. The author would also like to thank Professor Almut Burchard and Professor Larry Guth for helpful suggestions. \textcopyright 2011 by
the authors.}}

\author{Shibing Chen\thanks{Department of Mathematics, University of Toronto, Toronto, Ontario
Canada M5S 2E4 sbchen@math.toronto.edu.}}
\date{\today}

\begin{document}
 \maketitle

\begin{abstract}
We prove a one-parameter family of sharp integral inequalities for functions on the $n$-dimensional unit ball. The inequalities are conformally invariant, and the sharp constants are attained for functions that are equivalent to a constant function under conformal transformations. As a limiting case, we obtain an inequality that generalizes Carleman's inequality for harmonic functions in the plane to poly-harmonic functions in higher dimensions.
\end{abstract}

\section{\textbf{Introduction}}
In this paper we compute the norm and the maximizers of some integral operators between Lebesgue function spaces $L^p$, used to solve the Dirichlet problems for certain elliptic partial differential equations on the ball and the upper halfspace. One of the difficulties confronted is a lack of compactness due to conformal invariance of the equations in question. In some limiting cases, we are able to differentiate these inequalities with respect to the parameter $1/p$ to obtain nonlinear weighted exponential-type inequalities governing solutions and subsolutions of higher order linear elliptic problems on the ball, which improve on the maximum principle and give higher dimensional generalizations of Carleman's isoperimetric inequality for subharmonic functions on the unit disc.

In the following $B_{n}$ denotes the $n$-dimensional unit ball in Euclidean space,
$\|u\|_{L^{p}(\Omega)}$ is the $L^{p}$ norm of function $u$ defined on domain $\Omega$,
$|B_{n}|$ is the volume of $B_{n}$ and $c(n,a,p)$ is some constant which depends on $n,a$ and $p$. The parameter $a$ satisfies $2-n<a<1$. Before giving the main theorems, we will give an interesting corollary for the reason that it is easy to state and it is clearly a natural
generalization of Carleman's inequality (see (8) below) in four dimension.

\begin{corollary 1}
For any $u$ satisfying $\Delta^{2}u\leq 0$ and $-\frac{\partial u}{\partial\gamma}\leq 1$, where
$\gamma$ is the outer nomal of $\partial B_{4}$, we have
\begin{eqnarray}
(\int_{B_{4}}e^{4u}dx)^{\frac{1}{4}} \leq S(\int_{\partial B_{4}}e^{3u}d\xi)^{\frac{1}{3}},
\end{eqnarray}
The sharp constant is assumed by the solution of $\Delta^{2}u=0$ in $B_{n}$ with boundary values
$-\frac{\partial u}{\partial\gamma}=1$ and $u=0$ on $\partial B_{4}$.
\end{corollary 1}

For a function defined on $\textbf{R}^{n-1}$ (thought of as the boundary of the upper half-space $\textbf{R}^{n}_{+}),$
we define a poly-harmonic extension as follows: for $(X,x_{n})\in \textbf{R}^{n}_{+}=\textbf{R}^{n-1}\times (0,+\infty)$,
\begin{eqnarray}
(P_{a}f)(X,x_{n})=d_{n,a}\int_{\textbf{R}^{n-1}}\frac{x_{n}^{1-a}}{((X-Y)^{2}+x_{n}^{2})^{\frac{n-a}{2}}}f(Y)d^{n-1}Y.
\end{eqnarray}
Here the choice of integrand guarantees independence of $P_{a} 1$ on $(X,x_{n})$, while $a<1$ ensures $P_{a} 1<\infty$, and
the normalization constants $d_{n,a}$ are chosen so that $P_{a} 1 = 1$ (and can be expressed explicitly using $\Gamma$ functions).
Recalling that inversion in the unit sphere maps the halfspace $x_{n}>1/2$ to the unit ball centered at $(\textbf{0},1)$,
we see the conformal map
\begin{eqnarray}
\phi(X,x_{n}) = \frac{(X,x_{n}+\frac{1}{2})}{|(X,x_n+\frac{1}{2})|^2}- (\textbf{0},1)
\end{eqnarray}
 maps the upper halfspace $x_{n}>0$ to
the standard ball $\phi:\textbf{R}^{n}_{+} \longrightarrow B_{n}$. Conformality of this map makes it easy to
compute its Jacobian $J(\phi) = |(X, x_{n}+\frac{1}{2})|^{-2n}$, and the Jacobian $J(\phi|_{\partial \textbf{R}^{n}_{+}}) = |(X,\frac{1}{2})|^{-2(n-1)}$ of
its boundary trace. Indeed, $\phi$ pulls back the Euclidean metric $g$ on $B_{n}$ to the conformally flat
metric $\phi^{*} g = |(X,x_{n}+\frac{1}{2})|^{-4} \sum dx_i^2$ on $\textbf{R}^{n}_{+}$. Then it is not hard to check the formula
\begin{eqnarray}
f(X,x_{n}) = |(X,x_{n}+\textstyle\frac{1}{2})|^{2-n-a} \tilde f \circ \phi(X,x_{n})
\end{eqnarray}
 and its restriction to $x_n=0$ boundary trace define Banach space isometries from $\tilde f \in L^{\frac{2n}{n-2+a}}(B_{n})$ to $f \in L^{\frac{2n}{n-2+a}}(\textbf{R}^{n}_{+})$ and from $L^{\frac{2(n-1)}{n-2+a}}(\partial B_{n})$ to $L^{\frac{2(n-1)}{n-2+a}}(\textbf{R}^{n-1})$ respectively. We define the poly-harmonic extension $\tilde P_a \tilde f$ of $\tilde f \in L^{\frac{2(n-1)}{n-2+a}}(\partial B_n)$ implicitly by using $P_a$ after pulling back from the ball to the halfspace:
\begin{eqnarray}(\tilde P_a \tilde f) \circ \phi(X,x_n) = |(X,x_n+{\textstyle \frac{1}{2}})|^{n+a-2}P_a\left(\frac{\tilde f \circ \phi(Y,\frac{1}{2})}{|(Y,\frac{1}{2})|^{n+a-2}}\right).
\end{eqnarray}
 When $a=0$, $\tilde P_a \tilde f$ again becomes the usual harmonic extension to the ball. Another case of special interest is $a=2-n$, in which case the conformal factors are suppressed so that $\tilde P_{2-n} 1 = 1$, and the isometric Banach spaces are both of $L^\infty$ type. When $n=2k$ the extended function turns out to be $k$ harmonic on the $2k$ dimensional ball, i.e. $\Delta^k \tilde P_{2-2k} \tilde f =0$.

\begin{theorem}

For any $f\in L^{\frac{2(n-1)}{n-2+a}}(\partial B_{n})$, $n\geq 2$ and $n-2+a>0$, we have the sharp inequality
\begin{eqnarray}
\|\widetilde{P}_{a}f\|_{L^{\frac{2n}{n-2+a}}(B_{n})}\leq S_{n,a}\|f\|_{L^{\frac{2(n-1)}{n-2+a}}(\partial B_{n})},
\end{eqnarray}
 where the sharp constant $S_{n,a}$ depends only on $n$ and $a$. The optimizers are unique up to a conformal transform
 and include the constant function $f=1$.

\end{theorem}
We now study the limiting information. Letting $f=1+\frac{n-2+a}{2}F$ and $a\rightarrow 2-n$, we get the following inequality

\begin{theorem}

For any $F$ such that $e^{F}\in L^{n-1}(\partial B_{n})$, $n> 2$, we have
\begin{eqnarray}
\|e^{I_{n}+\widetilde{P}_{2-n}F}\|_{L^{n}( B_{n})}\leq S_{n} \|e^{F}\|_{L^{n-1}(\partial B_{n})},
\end{eqnarray}
where $I_{n}=\left(\log(X^{2}+(x_{n}+\frac{1}{2})^{2})-d_{n,2-n}\int_{R^{n-1}}\frac{x_{n}^{n-1}}{((X-Y)^{2}+x_{n}^{2})^{n-1}}\log(Y^{2}+\frac{1}{4})d^{n-1}Y\right)\circ \phi^{-1}.$ Up to a conformal transform any constant is an optimizer.

\end{theorem}

In Theorem 2, when $n=4$ we obtain Corollary 1.

\section{\textbf{Some related literature}}
Theorem 2 is a natural higher dimensional generalization of a well known inequality by Carleman [1]
\begin{eqnarray}
\int_{B_{2}}e^{2u}dx\leq\frac{1}{4\pi}(\int_{\partial B_{2}}e^{u}d\theta)^{2},
\end{eqnarray}
for all harmonic functions in $B_{2}$. Equality occurs exactly for $u=c$ and $u=-2\log|x-x_{0}|+c$, where $c$ is a constant and $x_{0}\in R^{2}-\overline{B_{2}}.$

Although Carleman proved (8) initially for harmonic functions, it follows from the maximum principle that inequality (8) holds for subharmonic functions.
Beckenbach and Rado [6] used Carleman's inequality to prove the isoperimetric inequality on a surface with non-positive Gauss curvature:
If on a surface with non-positive Gauss curvature an analytic curve $C$ of length $L$ encloses a simply-connected domain $D$ of area $A$, then the inequality
$$L^{2}\geq 4\pi A$$
holds. This is exactly the sharp isoperimetric inequality in the plane. Their proof is quite simple: In isothermal coordinates  $(x,y)$ for a simply-connected domain $\widetilde{D}$ which is slightly larger than $D$,
then the metric on $\widetilde{D}$ can be written as $e^{2w}(dx^{2}+dy^{2})$, for $(x,y)$ in some bounded domain $\Omega\in \textbf{R}^{2}$.
Now, the coordinate image of $D$ in $\Omega$  is a Jordan domain, so by the Riemann mapping theorem we can map it to $B_{2}$ conformally, which means
$D$ with the metric induced by the metric of the surface is isometric to $(B_{2},e^{2u}g)$, where $u$ is a subharmonic function (By the non-positive
curvature condition). Beckenbach and Rado's result now follows directly from Carleman's inequality.

The 2008 paper by Hang, Wang and Yan [2], generalized this inequality to higher dimensions as follows. For any harmonic function $u$ in $B_{n}$,
\begin{eqnarray}
\|u\|_{L^{\frac{2n}{n-2}}(B_{n})}\leq n^{-\frac{n-2}{2(n-1)}}\omega_{n}^{-\frac{n-2}{2n(n-1)}}\|u\|_{L^{\frac{2(n-1)}{n-2}}(\partial B_{n})},
\end{eqnarray}
where $n\geq 3$ and $\omega_{n}$ is the volume of $B_{n}$. Any constant is an optimizer and it is unique up to a conformal transformation (as will be explained before the proof of Theorem 2).
This is a special case of our Theorem 1 (a=0), and again because of the maximum principle, this inequality holds for subharmonic functions. Hang,
Wang and Yan interpreted their inequality as the isoperimetric inequality for $B_{n}$ with metric $\rho^{\frac{4}{n-2}}g$, where $\rho$ is subharmonic (which means non-positive scalar curvature).
By using the conformal map (2), the equivalent form of inequality (9) in the upper-half space is
\begin{eqnarray}
\|P_{0}f\|_{L^{\frac{2n}{n-2}(\textbf{R}_{+}^{n})}}\leq n^{-\frac{n-2}{2(n-1)}}\omega_{n}^{-\frac{n-2}{2n(n-1)}}\|f\|_{L^{\frac{2(n-1)}{n-2}}(\textbf{R}^{n-1} )},
\end{eqnarray}
for all $f\in L^{\frac{2(n-1)}{n-2}}(\textbf{R}^{n-1}).$
Here $\textbf{R}^{n-1}$ is the boundary of $\textbf{R}_{+}^{n}$. The optimizers are $f(Y)=\frac{c}{(\lambda^{2}+|Y-Y_{0}|^{2})^\frac{n-2}{2}}$ , for some constant c,
positive constant $\lambda$ and $Y_{0}\in \textbf{R}^{n-1}$.

\textbf{Remark 1.} We point out that the sharp inequality (9) combines with Brezis and Lieb's dual argument ([8] page10-11) to give the sharp
version of inequality (1.9) in [8] when the domain is a ball: $$\|\nabla f\|_{L^{2}(B_{n})}+  C(n)\|f\|_{L^{\frac{2(n-1)}{n-2}}(\partial B_{n})}
\geq S_{n}\|f\|_{L^{\frac{2n}{n-2}}(B_{n})},$$ where $S_{n}$ is sharp Sobolev constant and $C(n)$ can be determined by letting $f=1$ when the inequality becomes
equality. This sharp Sobolev inequality with trace term was also proved by Maggi and Villani in [9] by using methods from optimal transportation.

\textbf{Remark 2.} When $-1<a<1$, from Caffarelli and Silvestre [4] we know $u=P_{a}f$ is the unique solution to the boundary value problem
\begin{eqnarray*}
\mathrm{div}(x_{n}^{a}\nabla u(X,x_{n}))&=&0, (X,x_{n})\in \textbf{R}^{n}_{+}\\
u(X,0)&=&f, X\in \textbf{R}^{n-1}.
\end{eqnarray*}
Then the fractional Laplacian can be defined by using an analogue of the Dirichlet to Neumann map $(-\Delta)^{\frac{1-a}{2}}f=-\displaystyle{\lim_{x_{n}\rightarrow 0+}}x_{n}^{a}u_{y}$. So, our equivalent form of inequality (6) on $\textbf{R}^{n}_{+},$ namely
$\|P_{a}f\|_{L^{\frac{2n}{n-2+a}}(\textbf{R}^{n}_{+})}\leq \widetilde{S}_{n,a}\|f\|_{L^{\frac{2(n-1)}{n-2+a}}(\textbf{R}^{n})},$ provides a sharp estimate for the
$L^{\frac{2n}{n-2+a}}$ norm of  solution of the above boundary value problem.

\section{\textbf{Proof of Theorem 1}}

Since $P_{a}$ enjoys very similar properties to the special case $P_{0}$ (classical harmonic extension),
we are also able to use the method of symmetrization developed by Lieb [7] to prove the existence of maximizer as Hang, Wang and Yan did in [2].
The following Lemmas are parallel to those in [2], but notice that now we are dealing with poly-harmonic extension instead of harmonic extension.

Recall if $X$ is a measure space, $p > 0$ and $u$ is a measurable function on $\Omega$, then
$$\|u\|_{L^{p}_{w}}=\sup_{t>0} t||u|>t|^{\frac{1}{p}}.$$
The space $L^{p}_{w}(X)=\{$$u$:$u$ is measurable and $\|u\|_{L^{p}_{w}(\Omega)}<\infty\}$. More generally, for any $0<p<\infty$ and $0<q\leq\infty$,
we have Lorentz norm $\|\cdot\|_{L^{p,q}}$ which is defined by $\|u\|_{L^{p,q}}=p^{\frac{1}{q}}(\int_{0}^{\infty}t^{q}||u|\geq t|^{\frac{q}{p}}\frac{dt}{t})^{\frac{1}{q}}$
and Lorentz space $L^{p,q}(\Omega)$ . $L^{p}_{w}(\Omega)=L^{p,\infty}(\Omega)$ is a special case of such spaces.

\begin{lemma}Defining $P_{a}$ as in (1), there exist constants $c_{n,a}$ and $c_{n,a,p}$ such that
$$\|P_{a}f\|_{L^{\frac{n}{n-1}}_{w}(\textbf{R}^{n}_{+})}\leq c_{n,a}\|f\|_{L^{1}(\textbf{R}^{n-1})}$$
and
$$\|P_{a}f\|_{L^{\frac{np}{n-1}}(\textbf{R}^{n}_{+})}\leq c_{n,a,p}\|f\|_{L^{p}(\textbf{R}^{n-1})}$$
for all $1<p\leq\infty$.
Moreover for $1<p<\infty$ we have
$$\|P_{a}f\|_{L^{\frac{np}{n-1}}(\textbf{R}^{n}_{+})}\leq c_{n,a,p}\|f\|_{L^{p,\frac{np}{n-1}}(\textbf{R}^{n-1})}.$$
\end{lemma}
Proof of Lemma 1. To prove the weak
estimate, we may assume $f\geq0$ and $\|f\|_{L^{1}(\textbf{R}^{n-1})}=1$. It is easy to see
$(P_{a}f)(X,x_{n})\leq \frac{d_{n,a}}{x_{n}^{n-1}}$ for $(X,x_{n})\in \textbf{R}^{n}_{+}$ and
\begin{eqnarray*}
&&\int_{(X,x_{n})\in \textbf{R}^{n}_{+},0<y<b}(P_{a}f)(X,x_{n})d^{n-1}Xdx_{n}\\
&=&\int_{\textbf{R}^{n-1}}d^{n-1}Y(f(Y)\int_{0}^{b}dy\int_{\textbf{R}^{n-1}}d_{n,a}\frac{x_{n}^{1-a}}{((X-Y)^{2}+x_{n}^{2})^{\frac{n-a}{2}}}d^{n-1}X)\\
&=&b
\end{eqnarray*}
for $b>0$. Hence for $t>0$,
\begin{eqnarray*}
|P_{a}f>t|&=&|\{(X,x_{n})\in \textbf{R}^{n}_{+}:0<x_{n}<(d_{n,a}^{-1}t)^{-\frac{1}{n-1}},(P_{a}f)(X,x_{n})>t\}|\\
&\leq& \frac{1}{t}\int_{X\in \textbf{R}^{n-1},0<x_{n}<(d_{n,a}^{-1}t)^{-\frac{1}{n-1}}}(P_{a}f)(X,x_{n})d^{n-1}Xdx_{n}\\
&=&(d_{n,a}^{-1}t)^{-\frac{1}{n-1}}
\end{eqnarray*}
The weak type inequality follows. The strong estimate follows
from Marcinkiewicz interpolation theorem (see [5], p197) and the basic fact
$\|P_{a}f\|_{L^{\infty}(\textbf{R}^{n}_{+})}\leq \|f\|_{L^{\infty}(\textbf{R}^{n-1})}$.

\textbf{Remark 3.}  In fact when $p=\frac{2(n-1)}{n-2}$ and $a=0$, the second estimate was also proved by Brezis and Lieb [8] by using some elementary dual argument.

\begin{lemma}
If $n\geq 2$ and $1<p<\infty$, then the supremum
\begin{eqnarray}
c_{n,a,p}^{\frac{np}{n-1}}=\sup\{\|P_{a}f\|_{L^{\frac{np}{n-1}}}^{\frac{np}{n-1}}: \|f\|_{L^{p}(\textbf{R}^{n-1})}=1\},
\end{eqnarray}
is attained by some function. After multiplying by a nonzero constant, every maximizer $f$ is nonnegative,
radially symmetric with respect to some point, strictly decreasing
in the radial direction and it satisfies the following Euler-Lagrange equation
\begin{eqnarray*}
f(Y)^{p-1}=\int_{\textbf{R}^{n}_{+}}\frac{x_{n}^{1-a}}{((X-Y)^{2}+x_{n}^{2})^{\frac{n-a}{2}}}(P_{a}f)^{\frac{np}{n-1}-1}(X,x_{n})d^{n-1}X dx_{n}.
\end{eqnarray*}
In particular, if $n\geq 2$, $p=\frac{2(n-1)}{n-2+a}$ and $n-2+a>0$, then every maximizer is of the form
\begin{eqnarray}
f(Y)=\pm c(n,a)(\frac{\lambda}{\lambda^{2}+|Y-Y_{0}|^{2}})^{\frac{n-2+a}{2}}
\end{eqnarray}
for some $\lambda>0$,
$Y_{0}\in \textbf{R}^{n-1}$.
\end{lemma}%
\emph{Proof of Lemma 2.} %
First we recall the important Riesz rearrangement inequality.
Let $u$ be a measurable function on $\textbf{R}^{n}$, the symmetric rearrangement of $u$ is
 the nonnegative lower semi-continuous radial decreasing function $u^{*}$ that has the
same distribution as $u$. We have
$$\int_{\textbf{R}^{n}}dx\int_{\textbf{R}^{n}}u(x)v(y-x)w(y)dy\leq \int_{\textbf{R}^{n}}dx\int_{\textbf{R}^{n}}u^{*}(x)v^{*}(y-x)w^{*}(y)dy.$$
Using the fact $\|w\|_{L^{p}(\textbf{R}^{n})}=\|w^{*}\|_{L^{p}(\textbf{R}^{n})}$ for $p>0$ and the standard duality argument, we see for $1\leq p\leq \infty,$
$$\|u*v\|_{L^{p}(\textbf{R}^{n})}\leq \|u^{*}*v^{*}\|_{L^{p}(\textbf{R}^{n})}.$$
Moreover if $u$ is nonnegative radially symmetric and strictly decreasing in the radial
direction, $v$ is nonnegative, $1<p<\infty$ and
$$\|u*v\|_{L^{p}(\textbf{R}^{n})}= \|u^{*}*v^{*}\|_{L^{p}(\textbf{R}^{n})}<\infty,$$
then for some $x_{0}\in \textbf{R}^{n},$  we have $v(x)=v^{*}(x-x_{0})$.

Now, assume $f_{i}$ is a maximizing sequence in (11). Since $\|f_{i}^{*}\|_{L^{p}(\textbf{R}^{n-1})}=\|f_{i}\|_{L^{p}(\textbf{R}^{n-1})}=1$ and
\begin{eqnarray*}
\|P_{a}f_{i}\|_{L^{\frac{np}{n-1}}(\textbf{R}^{n}_{+})}^{\frac{np}{n-1}}&=&\int_{0}^{\infty}\|P_{a,x_{n}}*f_{i}\|_{L^{\frac{np}{n-1}}(\textbf{R}^{n-1})}^{\frac{np}{n-1}}dy\\
&\leq &\int_{0}^{\infty}\|P_{a,x_{n}}*f^{*}_{i}\|_{L^{\frac{np}{n-1}}(\textbf{R}^{n-1})}^{\frac{np}{n-1}}dx_{n}\\
&=&\|P_{a}f^{*}_{i}\|_{L^{\frac{np}{n-1}}(\textbf{R}^{n}_{+})}^{\frac{np}{n-1}},
\end{eqnarray*}
where $P_{a,x_{n}}=d_{n,a}\frac{x_{n}^{1-a}}{(X^{2}+x_{n}^{2})}$ and notice that it is symmetric and strictly decreasing in the radial
direction of $X$ variable for any fixed $x_{n}$. We see $f_{i}$ is again a maximizing sequence. Hence we may assume $f_{i}$ is a nonnegative
radial decreasing function.

For any $f\in L^{p}(\textbf{R}^{n-1})$ and any $\lambda >0$, we let $f^{\lambda}(Y)=\lambda^{-\frac{n-1}{p}}f(\frac{Y}{\lambda})$, so that is clear that
$(P_{a}f^{\lambda})(X,x_{n})=\lambda^{-\frac{n-1}{p}}(P_{a}f)(\frac{X}{\lambda},\frac{x_{n}}{\lambda})$ and hence
 $\|f^{\lambda}\|_{L^{p}(\textbf{R}^{n-1})}=\|f\|_{L^{p}(\textbf{R}^{n-1})}$ and
 $\|P_{a}f^{\lambda}\|_{L^{\frac{np}{n-1}}(\textbf{R}^{n}_{+})}=\|P_{a}f\|_{L^{\frac{np}{n-1}}(\textbf{R}^{n}_{+})}$. For convenience, denote $e_{1}=(1,0,\cdots,0)\in \textbf{R}^{n-1}$
and
$$a_{i}=\sup\{f_{i}^{\lambda}(e_{1})|\lambda>0\}=\sup\{\lambda^{-\frac{n-1}{p}}f_{i}(\frac{e_{1}}{\lambda})|\lambda>0\}.$$
It follows that $0\leq f_{i}(Y)\leq a_{i}|Y|^{-\frac{n-1}{p}},$
and hence $\|f_{i}\|_{L^{p,\infty}(\textbf{R}^{n-1})}\leq |B_{n-1}|^{\frac{1}{p}}a_{i}$.

Now
\begin{eqnarray*}
  \|P_{a}f_{i}\|_{L^{\frac{np}{n-1}}(\textbf{R}^{n}_{+})}&\leq& c(n,a,p)\|f_{i}\|_{L^{p,\frac{np}{n-1}}(\textbf{R}^{n-1})}\\
  &\leq& c(n,a,p)\|f_{i}\|_{L^{p}(\textbf{R}^{n-1})}^{\frac{n-1}{n}}\|f_{i}\|^{\frac{1}{n}}_{L^{p,\infty}(\textbf{R}^{n-1})}\\
  &\leq& c(n,a,p)a_{i}^{\frac{1}{n}},
\end{eqnarray*}
which implies $a_{i}\geq c(n,a,p)>0$. We may choose $\lambda_{i}>0$ such that $f_{i}^{\lambda_{i}}(e_{1})\geq c(n,a,p)>0$.
Replacing $f_{i}$ by $f_{i}^{\lambda_{i}}$ we may assume $f(e_{1})\geq c(n,a,p)>0$. On the other hand , since $f_{i}$ is nonnegative radial
decreasing and $\|f_{i}\|_{L^{p}(\textbf{R}^{n-1})}=1$, we see
$$|f_{i}(Y)|\leq|B_{n-1}|^{-\frac{1}{p}}|Y|^{-\frac{n-1}{p}}.$$
Hence after passing to a subsequence, we may find a nonnegative radial decreasing function $f$ such that $f_{i} \rightarrow f $ a.e. It follows
that $f(Y)\geq c(n,a,p)>0$ for $|Y|\leq 1$, $f_{i}\rightharpoonup f$ in $L^{p}(\textbf{R}^{n-1})$ and $\|f\|_{L^{p}(\textbf{R}^{n-1})}\leq 1$. By
Lieb [7](Lemma 2.6), we have
$$ \int_{\textbf{R}^{n-1}}||f_{i}|^{p}-|f|^{p}-|f_{i}-f|^{p}|d^{n-1}Y\rightarrow0.$$
It follows that
\begin{eqnarray*}\|f_{i}-f\|^{p}_{L^{p}(\textbf{R}^{n-1})}&=&\|f_{i}\|^{p}_{L^{p}(\textbf{R}^{n-1})}-\|f\|^{p}_{L^{p}(\textbf{R}^{n-1})}+o(1)\\
&=&1-\|f\|^{p}_{L^{p}(\textbf{R}^{n-1})}+o(1).
\end{eqnarray*}
On the other hand, since $(P_{a}f_{i})(X,x_{n})\rightarrow (P_{a}f)(X,x_{n}))$ for $(X,x_{n})\in \textbf{R}^{n}_{+}$
and $\|P_{a}f_{i}\|_{L^{\frac{np}{n-1}}(\textbf{R}^{n}_{+})}\leq c_{n,a,p}$, we see
\begin{eqnarray*}\|P_{a}f_{i}\|^{\frac{np}{n-1}}_{L^{\frac{np}{n-1}}(\textbf{R}^{n}_{+})}&=&\|P_{a}f\|^{\frac{np}{n-1}}_{L^{\frac{np}{n-1}}(\textbf{R}^{n}_{+})}
+ \|P_{a}f_{i}-P_{a}f\|^{\frac{np}{n-1}}_{L^{\frac{np}{n-1}}(\textbf{R}^{n}_{+})}+o(1)\\
&\leq& c_{n,a,p}^{\frac{np}{n-1}}\|f\|^{\frac{np}{n-1}}_{L^{p}(\textbf{R}^{n-1})}+c_{n,a,p}^{\frac{np}{n-1}}\|f_{i}-f\|^{\frac{np}{n-1}}_{L^{p}(\textbf{R}^{n-1})}+o(1).\\
\end{eqnarray*}
Hence
$$1\leq \|f\|^{\frac{np}{n-1}}_{L^{p}(\textbf{R}^{n-1})}+\|f_{i}-f\|^{\frac{np}{n-1}}_{L^{p}(\textbf{R}^{n-1})}+o(1).$$
Let $i\rightarrow\infty$, we see
$$1\leq \|f\|^{\frac{np}{n-1}}_{L^{p}(\textbf{R}^{n-1})}+(1-\|f\|^{p}_{L^{p}(\textbf{R}^{n-1})})^{\frac{n}{n-1}}.$$
Since$\frac{n}{n-1}>1$ and $f\neq 0$, we see $\|f\|_{L^{p}(\textbf{R}^{n-1})}=1$. Hence $f_{i}\rightarrow f$ in $L^{p}(\textbf{R}^{n-1})$ and $f$ is
a maximizer. This implies the existence of an extremal function.

Assume $f\in L^{p}(\textbf{R}^{n-1})$ is a maximizer, then so is $|f|$.
Hence $\|P_{a}f\|_{L^{\frac{np}{n-1}}(\textbf{R}^{n}_{+})}=\|P_{a}|f|\|_{L^{\frac{np}{n-1}}(\textbf{R}^{n}_{+})}$. On the other hand ,
since $|(P_{a}f)(X,x_{n})|\leq (P_{a}|f|)(X,x_{n})$ for $(X,x_{n})\in \textbf{R}^{n}_{+}$, we see $|P_{a}f|=P_{a}(|f|)$ and this implies either $f\geq 0$ or
$f\leq 0$. Assume $f\geq 0$, then the Euler-Lagrange equation after scaling by a positive constant is given by
$$f(Y)^{p-1}=\int_{\textbf{R}^{n}_{+}}\frac{y^{1-a}}{((X-Y)^{2}+x_{n}^{2})^{\frac{n-a}{2}}}(P_{a}f)^{\frac{np}{n-1}-1}(X,x_{n})d^{n-1}Xdx_{n}.$$
On the other hand, we know for $x_{n}>0$,
$$\|P_{a,x_{n}}*f\|_{L^{\frac{np}{n-1}}(\textbf{R}^{n-1})}=\|P_{a,x_{n}}*f^{*}\|_{L^{\frac{np}{n-1}}(\textbf{R}^{n-1})}$$
which implies $f(Y)=f^{*}(Y-Y_{0})$ for
some $Y_{0}$. It follows from the above Euler-Lagrange equation and Lemma 2.2 of Lieb [7] that $f$ must be strictly decreasing along the radial direction.

For the case when $p=\frac{2(n-1)}{n-2+a}$, we first observe that if $f\in L^{\frac{2(n-1)}{n-2+a}}(\textbf{R}^{n-1})$, let $u=P_{a}f$,
$\widetilde{f}=\frac{1}{|Y|^{n-2+a}}f(\frac{Y}{|Y|^{2}})$ and $\widetilde{u}=\frac{1}{|(X, x_{n})|^{n-2+a}}f(\frac{(X, x_{n})}{|(X, x_{n})|^{2}})$,
then we have $\widetilde{u}=P_{a}\widetilde{f}$, $\|\widetilde{f}\|_{L^{\frac{2(n-1)}{n-2+a}}(\textbf{R}^{n-1})}=\|f\|_{L^{\frac{2(n-1)}{n-2+a}}(\textbf{R}^{n-1})}$
and $\|\widetilde{u}\|_{L^{\frac{2n}{n-2+a}}(\textbf{R}^{n}_{+})}=\|u\|_{L^{\frac{2n}{n-2+a}}(\textbf{R}^{n}_{+})}$ .
This is the conformal invariance property for the particular power. As a consequence, if $f$ is a maximizer which is nonnegative and radial,
then$ \frac{1}{|Y|^{n-2+a}}f(\frac{Y}{|Y|^{2}}-e_{1})$ is also a maximizer. In particular, $ \frac{1}{|Y|^{n-2+a}}f(\frac{Y}{|Y|^{2}}-e_{1})$
is radial with respect to some points. To find such $f$, we need the following useful Proposition of Hang, Wang and Yan [2]( Proposition 4.1).

\begin{lemma}
Let $n\geq 2$, $u$ be a function on $\textbf{R}^{n}$ which is radial with respect to the origin, $0<u(x)<\infty$ for $x\neq 0,$ $e_{1}=(1,0,\cdots,0),$ $\alpha\in R,$
$\alpha\neq 0$. If $v(x)=|x|^{\alpha}u(\frac{x}{|x|^{2}}-e_{1})$ is radial with respect to some point, then either $u(x)=(c_{1}|x|^{2}+c_{2})^{\frac{\alpha}{2}}$
for some $c_{1}\geq 0, c_{2}>0$ or
\begin{equation}
u(x)=\begin{cases}
c_{1}|x|^{\alpha} & \text{if}\ x\neq 0 \\
c_{2} & \text{if}\ x=0.
\end{cases}
\end{equation}
\end{lemma}

\emph{Proof of Lemma 2 continued.} Since $\|f\|_{L^{\frac{2(n-1)}{n-2+a}}(\textbf{R}^{n-1})}=1$ and it is strictly decreasing
along the radial direction, we have $0<f(Y)<\infty$ for all $Y\neq 0$. Note that since $f$ satisfies
the Euler-Lagrange equation, it is defined everywhere instead of almost everywhere. It follows from Lemma 3 that
$f(Y)=(c_{1}|Y|^{2}+c_{2})^{-\frac{n-2+a}{2}}$ for some $c_{1}, c_{2}>0$ (since $f$ can not be constant function and the scalar multiple
of $|Y|^{2-n}$ is ruled out by the integrability). Using the condition $\|f\|_{L^{\frac{2(n-1)}{n-2+a}}(\textbf{R}^{n-1})}=1$, it is easy to see
$c_{1}c_{2}=c_{n,a}$. Hence for some $\lambda>0$,
$$f(Y)=c(n,a)(\frac{\lambda}{\lambda^{2}+|Y-Y_{0}|^{2}})^{\frac{n-2+a}{2}}.$$

\emph{Proof of Theorem 1.} For any $\widetilde{f}\in L^{\frac{2(n-1)}{n-2+a}}(\partial B_{n}),$
let $\widetilde{u}=\widetilde{P}_{a}f, $
$$f=\frac{1}{|(Y,0)+(\textbf{0},\frac{1}{2})|^{n-2+a}}\widetilde{f}\circ\phi,$$
 and
$$u=\frac{1}{|(X,x_{n})+(\textbf{0},\frac{1}{2})|^{n-2+a}}\widetilde{u}\circ\phi.$$
By definition (4) we have $u=P_{a}f$ and by the discussion below (3)
we have $\|\widetilde{f}\|_{L^{\frac{2(n-1)}{n-2+a}}(\partial B_{n})}=\|f\|_{L^{\frac{2(n-1)}{n-2+a}}(\textbf{R}^{n-1})}$ and
$\|\widetilde{u}\|_{L^{\frac{2n}{n-2+a}}( B_{n})}=\|u\|_{L^{\frac{2n}{n-2+a}}( \textbf{R}^{n}_{+})}$. Then, Theorem 1 follows easily from the above facts and
Lemma 2.

\section{\textbf{Proof of Theorem 2}}
First we will discuss some conformal invariance properties of the operator $\widetilde{P}_{a}$.
Let $\widetilde{\tau}$ be a conformal transform from $B_{n}$ to itself, $\tau=\widetilde{\tau}|_{\partial B_{n}}$ is the induced conformal transform from
$\partial B_{n}$ to itself, $\widetilde{J}$ is the Jacobian of $\widetilde{\tau}$, $J$ is the Jacobian of $\tau$, $\varepsilon=n-2+a$.
For $f\in L^{\frac{2(n-1)}{n-2+a}}(\partial B_{n})$, when $\varepsilon\neq 0$, we have
\begin{eqnarray}
\widetilde{P}_{a}(J^{\frac{\varepsilon}{2(n-1)}}f\circ \tau)=\widetilde{J}^{\frac{\varepsilon}{2n}}(\widetilde{P}_{a}f)\circ\widetilde{\tau}.
\end{eqnarray}
It is straightforward to check this property by using the definition of $\widetilde{P}_{a}$ in (5).

Now, for smooth function $f$, when $\varepsilon$ goes to 0 it is obvious that
\begin{eqnarray}
\widetilde{P}_{2-n}(f\circ\tau)=(\widetilde{P}_{2-n}f)\circ\widetilde{\tau}.
\end{eqnarray}
By letting $f=1$ and taking derivative with respect to $\varepsilon$ at 0, we have
\begin{eqnarray}
\frac{d(\widetilde{P}_{a}1)}{d\varepsilon}|_{\varepsilon=0}+\widetilde{P}_{2-n}(\frac{1}{2(n-1)}\log J)=\frac{d(\widetilde{P}_{a}1)}{d\varepsilon}|_{\varepsilon=0}\circ \widetilde{\tau}+\frac{1}{2n}\log\widetilde{J}
\end{eqnarray}
So the inequality in the Theorem 2 is invariant when $F$ is replaced by $F\circ\tau+\frac{1}{n-1}\log J$.

\emph{Proof of Theorem 2.}Recalling $\tilde P_{2-n} 1 = 1$, let $f=1+\varepsilon F$, where $F$ is some smooth function defined on $\partial B_{n}$.
By Theorem 1, we have the inequality
\begin{eqnarray*}
\|\widetilde{P}_{a}(1+\varepsilon F)\|_{L^{\frac{2n}{n-2+a}}(B_{n})}\leq S_{n,a}\|1+\varepsilon F\|_{L^{\frac{2(n-1)}{n-2+a}}(\partial B_{n})},
\end{eqnarray*}
which means
\begin{eqnarray*}
(\int_{B_{n}}(\widetilde{P}_{a}1)^{\frac{2n}{\epsilon}}(1+\frac{\varepsilon\widetilde{P}_{a}F}{\widetilde{P}_{a}1})^{\frac{2n}{\varepsilon}}dx)^{\frac{1}{n}}\leq S_{n,a}^{\frac{2}{\epsilon}}(\int_{\partial B_{n}}(1+\varepsilon F)^{\frac{2(n-1)}{\varepsilon}}d\xi)^{\frac{1}{n-1}}.
\end{eqnarray*}
Note that when $F=0$ the above inequality becomes equality, then by the following estimates we will see in this case the integrals in both sides will converge to some finite numbers, which means the constant $S_{n,a}^{\frac{2}{\epsilon}}$ will also converge.

In order to take limit $\varepsilon\rightarrow 0$, we need to apply the Dominated Convergence Theorem. we will bound the term $\widetilde{P}_{a}1$ from below by a constant $A$ and bound $(\widetilde{P}_{a}1)^{\frac{2n}{\varepsilon}}$ from above by a constant $B$, both $A$ and $B$ are  independent of $\varepsilon$. Let us derive the lower bound of $\widetilde{P}_{a}1$ first.
From (1) and (4) we know
\begin{eqnarray*}
(\widetilde{P}_{a}1)\circ\phi=(X^{2}+(x_{n}+\frac{1}{2})^{2})^{\frac{\varepsilon}{2}}
d_{n,a}\int_{\textbf{R}^{n-1}}\frac{x_{n}^{1-a}}{((X-Y)^{2}+x_{n}^{2})^{\frac{n-a}{2}}}\frac{1}
{(Y^{2}+\frac{1}{4})^{\frac{\varepsilon}{2}}}d^{n-1}Y,
\end{eqnarray*}
by letting $U=\frac{Y-X}{x_{n}}$ in the integral, we have
\begin{eqnarray*}
(\widetilde{P}_{a}1)\circ\phi&=&d_{n,a}\int_{\textbf{R}^{n-1}}\frac{1}{(u^{2}+1)^{\frac{n-a}{2}}}
\frac{(|X|^{2}+(x_{n}+\frac{1}{2})^{2})^{\frac{\varepsilon}{2}}}{((Ux_{n}+X)^{2}+\frac{1}{4})^{\frac{\varepsilon}{2}}}d^{n-1}U\\
&\geq&
d_{n,0}\int_{|u|\leq 1}\frac{1}{(|U|^{2}+1)^{n-1}}\frac{(|X|^{2}+(x_{n}+\frac{1}{2})^{2})^{\frac{\varepsilon}{2}}}{(2|X|^{2}+2x_{n}^{2}+\frac{1}{4})^{\frac{\varepsilon}{2}}}d^{n-1}U\\
&\geq&
d_{n,0}\int_{|U|\leq 1}\frac{1}{(|U|^{2}+1)^{n-1}}\frac{1}{2}d^{n-1}U\\
&=& A.
\end{eqnarray*}
Then, lets derive the upper bound of $(\widetilde{P}_{a}1)^{\frac{2n}{\varepsilon}}$, it is enough to
prove that $(\widetilde{P}_{a}1)^{\frac{n-2}{\varepsilon}}$ is bounded from above by some constant $B$ independent of $\varepsilon$.
As in the proof of lower bound, after the same change of variable we have
\begin{eqnarray*}
(\widetilde{P}_{a}1)^{\frac{n-2}{\varepsilon}}\circ\phi &=&
g(X,x_{n})\left(d_{n,a}\int_{\textbf{R}^{n-1}}\frac{1}{(|U|^{2}+1)^{\frac{n-a}{2}}}
\frac{1}{((Ux_{n}+X)^{2}+\frac{1}{4})^{\frac{\varepsilon}{2}}}d^{n-1}U\right)^{\frac{n-2}{\varepsilon}}\\
&\leq&
g(X,x_{n})d_{n,a}\int_{\textbf{R}^{n-1}}\frac{1}{(|U|^{2}+1)^{\frac{n-a}{2}}}
\frac{1}{((Ux_{n}+X)^{2}+\frac{1}{4})^{\frac{n-2}{2}}}d^{n-1}U\\
&\leq&
g(X,x_{n})\frac{d_{n,2-n}}{d_{n,0}}
d_{n,0}\int_{\textbf{R}^{n-1}}\frac{1}{(|U|^{2}+1)^{\frac{n}{2}}}
\frac{1}{((Ux_{n}+X)^{2}+\frac{1}{4})^{\frac{n-2}{2}}}d^{n-1}U\\
&=&
g(X,x_{n})\frac{d_{n,2-n}}{d_{n,0}}
d_{n,0}\int_{\textbf{R}^{n-1}}\frac{x_{n}}{((X-Y)^{2}+x_{n}^{2})^{\frac{n}{2}}}\frac{1}
{(|Y|^{2}+\frac{1}{4})^{\frac{n-2}{2}}}d^{n-1}Y\\
&=&
\frac{d_{n,2-n}}{d_{n,0}}\\
&=&
B,
\end{eqnarray*}
where $g(X,x_{n})=\left(|X|^{2}+(x_{n}+\frac{1}{2})^{2}\right)^{\frac{n-2}{2}}$. For the first inequality we applied Jensen's inequality, since $d_{n,a}\frac{1}{(u^{2}+1)^{\frac{n-a}{2}}}$
is a probability density in $\textbf{R}^{n-1}$ and $g(t)=t^{\frac{n-2}{\varepsilon}}$ is convex when $t\geq 0$.
The last identity holds because $$d_{n,0}\int_{\textbf{R}^{n-1}}\frac{x_{n}}{((X-Y)^{2}+x_{n}^{2})^{\frac{n}{2}}}\frac{1}
{(Y^{2}+\frac{1}{4})^{\frac{n-2}{2}}}d^{n-1}Y$$ is the harmonic extension of function $(Y^{2}+\frac{1}{4})^{-\frac{n-2}{2}}$ which is easy to verify that it is exactly $\frac{1}{(X^{2}+(x_{n}+\frac{1}{2})^{2})^{\frac{n-2}{2}}}$.

Now we can take limit $\varepsilon\rightarrow 0$ safely. By denoting
\begin{eqnarray*}
I_{n}&=& 2\frac{d(\widetilde{P}_{a}1)}{d\varepsilon}|_{\varepsilon=0}\\
&=&
\left(\log(X^{2}+(x_{n}+\frac{1}{2})^{2})-
d_{n,2-n}\int_{\textbf{R}^{n-1}}\frac{x_{n}^{n-1}}{((X-Y)^{2}+x_{n}^{2})^{n-1}}\log(Y^{2}+\frac{1}{4})d^{n-1}Y\right)\circ \phi^{-1},
\end{eqnarray*}
we get $\|e^{I_{n}+2\widetilde{P}_{2-n}F}\|_{L^{n}( B_{n})}\leq S_{n} \|e^{2F}\|_{L^{n-1}(\partial B_{n})}$.
After replacing $2F$ with $F$, the inequality in Theorem 2 is proved.
Since constant functions are optimizers for the above inequality, conformal invariance of the inequality  tells us that
the functions
$$F=C+\frac{1}{n-1}\log J$$
are also optimizers.

\textbf{Remark 4.} The uniqueness is lost when taking limit in the proof of Theorem 2. It would be  interesting to find a suitable method to prove that the optimizer is unique up to conformal transforms. The main difficulty seems to the author is that the integral kernel is too complicated.

\section{\textbf{Proof of Corollary 1}}
Now we are in the situation where $n=4$ and $a=-2$.
By [12] we know that $(P_{-2}f)(X,x_{n})=d_{4,-2}\int_{R^{3}}\frac{x_{n}^{3}}{((X-Y)^{2}+x_{n}^{2})^{n-1}}f(Y)dY$ is the bi-harmonic extension
of the function $f(Y)$ with boundary condition $\frac{\partial (P_{-2}f)(X,x_{n})}{\partial x_{n}}|_{x_{n}=0}=0$.
 It is straightforward to check that under the conformal map $\phi$ the bi-harmonic property and the Neumann boundary condition are preserved in dimension four, we have that
$\widetilde{P}_{-2}g$ is a bi-harmonic extension of a function $g$ defined on $S^{3}$ to a function on $B_{4}$ with boundary condition
$\frac{\partial \widetilde{P}_{-2}g}{\partial \gamma}|_{y=0}=0$. In view of Theorem 2, in order to prove Corollary 1 we only need to verify
$I_{2}$ satisfies $\Delta^{2}I_{2}= 0$, $I_{2}=0|_{S^{3}}$ and $-\frac{\partial I_{2}}{\partial\gamma}=1$.

From the formula for $I_{n}$, we have
\begin{eqnarray*}
I_{2}=\left(\log(|X|^{2}+(x_{n}+\frac{1}{2})^{2})-
d_{4,-2}\int_{R^{3}}\frac{x_{n}^{3}}{((X-Y)^{2}+x_{n}^{2})^{3}}\log(Y^{2}+\frac{1}{4})d^{n-1}Y\right)\circ \phi^{-1}.
\end{eqnarray*}
By using the explicit formula of $\phi$ one can get
\begin{eqnarray*}
I_{2}=2\log|\eta-S|-2C\int_{S^{3}}\frac{(1-|\eta|^{2})^{3}}{|\eta-\xi|^{6}}\log|\xi-S|d\xi,
\end{eqnarray*}
where $\eta$ is a point in $B_{4}$, $\xi$ is a point on $S^{3}$, $S$ is the south pole of $S^{3}$ and $C$ is the normalizing constant such that $C\int_{S^{3}}\frac{(1-|\eta|^{2})^{3}}{|\eta-\xi|^{6}}d\xi=1$.
Now from [12] we have the representation formula for bi-harmonic functions, namely for a smooth bi-harmonic function $g$ on $\overline{B^{4}}$ we have
\begin{eqnarray}
g(u)=C\int_{S^{3}}\frac{(1-|\eta|^{2})^{3}}{|\eta-\xi|^{6}}g(\xi)d\xi+
D\int_{S^{3}}\frac{(1-|\eta|^{2})^{2}}{|\eta-\xi|^{4}}(-\frac{\partial g}{\partial \gamma}(\xi)) d\xi,
\end{eqnarray}
where $D$ is a known constant.
Although the function $\log|\eta-S|$ is singular at the south pole $S$, If we apply the forthcoming approximation process we obtain
\begin{eqnarray*}
I_{2}&=&2\log|\eta-S|-2C\int_{S^{3}}\frac{(1-|\eta|^{2})^{3}}{|\eta-\xi|^{6}}\log|\xi-S|d\xi\\
&=&
-D\int_{S^{3}}\frac{(1-|\eta|^{2})^{2}}{|\eta-\xi|^{4}} d\xi,
\end{eqnarray*}
since $-\frac{\partial \log|\xi-S|}{\partial \gamma}(\xi)=-1$.
In the above equality $-D\int_{S^{3}}\frac{(1-|\eta|^{2})^{2}}{|\eta-\xi|^{4}} d\xi$ is the bi-harmonic extension of constant function 0 with boundary condition
$-\frac{\partial g}{\partial \gamma}(\xi)=1$, so $I_{2}$ satisfies all three conditions mentioned above.

Since $\log|\eta-S|$ is singular, we use approximation to justify the previous formula for $I_{2}$. Take a sequence $$S_{t}=(0,0,0,-t)\rightarrow S=(0,0,0,-1),$$ as $t\rightarrow 1+.$
Then $\log|\eta-S_{t}|$ is a smooth bi-harmonic function on $\overline{B^{4}}$, so we have
\begin{eqnarray*}
\log(\eta-S_{t})=C\int_{S^{3}}\frac{(1-|\eta|^{2})^{3}}{|\eta-\xi|^{6}}\log|\xi-S_{t}|d\xi
+D\int_{S^{3}}\frac{(1-|\eta|^{2})^{2}}{|\eta-\xi|^{4}}\frac{1+t y}{1+t^{2}+2ty}d\xi,
\end{eqnarray*}
here we use $y$ to denote the last coordinate of $\xi$.
For fixed $\eta\in B_{4}$, when $t$ approximates 1 from the right, $|\log|\xi-S_{t}||\leq |\log|\xi-S||$ for $\xi$ in a small neighborhood of $S$, since $|\log|\xi-S||$ is integrable on $S^{3}$, by the Dominated Convergence Theorem
$$C\int_{S^{3}}\frac{(1-|\eta|^{2})^{3}}{|\eta-\xi|^{6}}\log|\xi-S_{t}|d\xi\rightarrow C\int_{S^{3}}\frac{(1-|\eta|^{2})^{3}}{|\eta-\xi|^{6}}\log|\xi-S|d\xi,$$ as $t\rightarrow 1+$. Similarly, when $t$ is close to 1 from right hand side, we have $$|\frac{1+t y}{1+t^{2}+2ty}|\leq \frac{1}{2}+\frac{10}{|\xi-S|^{2}}.$$ Since
$\frac{1}{2}+\frac{10}{|\xi-S|^{2}}$ is integrable on $S^{3}$, by the Dominated Convergence Theorem again we have $$D\int_{S^{3}}\frac{(1-|\eta|^{2})^{2}}{|\eta-\xi|^{4}}\frac{1+t y}{1+t^{2}+2ty}d\xi\rightarrow \frac{1}{2}D\int_{S^{3}}\frac{(1-|\eta|^{2})^{2}}{|\eta-\xi|^{4}}d\xi,$$
as $t\rightarrow 1+$. Now by taking limit $t\rightarrow 1+$, it is clear that we have the representation formula for $\log|\eta-S|$. Finally since
the kernels in the representation formula (17) are positive, we conclude that the inequality in Corollary 1 is true for sub-biharmonic function $u$
with boundary conditions $-\frac{\partial u}{\partial\gamma}=1$ and $u=0$ on $\partial B_{4}$.

\end{document}